\documentclass[a4paper,twocolumn]{article} 

\usepackage[utf8]{inputenc}
\usepackage{amssymb,amsmath,amsthm} 
\usepackage{caption, subcaption} 
\usepackage{epsfig} 
\usepackage[numbers]{natbib} 
\usepackage{url} 
\newtheorem{prop}{Proposition} 
\newtheorem{lem}[prop]{Lemma} 
\newtheorem{cor}[prop]{Corollary} 
\usepackage[norelsize,algo2e,linesnumbered,ruled,vlined]{algorithm2e} 

\title{A dynamic programming approach to solving constrained linear-quadratic optimal control problems}

\author{Ruth Mitze and Martin M\"onnigmann\\
	Automatic Control and Systems Theory, Dep. of Mechanical Engineering,\\
	Ruhr-Universit\"at Bochum, 44801 Bochum, Germany.\\ E-mail: {\tt\small ruth.mitze@rub.de} and {\tt\small martin.moennigmann@rub.de}}

\begin{document}
\maketitle
          
\begin{abstract}                          
The solution of a constrained linear-quadratic regulator problem is determined by the set of its optimal active sets. We propose an algorithm that constructs this set of active sets for a desired horizon $N$ from that for horizon $N-1$. While it is not obvious how to extend the optimal feedback law itself for horizon $N-1$ to horizon $N$, a simple relation between the optimal active sets for two successive horizon lengths exists. Specifically, every optimal active set for horizon $N$ is a superset of an optimal active set for horizon $N-1$ if the constraints are ordered stage by stage. The stagewise treatment results in a favorable computational effort. In addition, it is easy to detect the solution of the current horizon is equal to the infinite-horizon solution, if such a finite horizon exists, with the proposed algorithm.  \\

\end{abstract}

\section{Introduction}\label{sec:Intro}
The constrained linear quadratic regulator (LQR) is solved by a piecewise-affine feedback law~\cite{Bemporad2002,Seron2003}.
This solution is conceptually simple and its piecewise-affine structure is the generalization of the optimal linear feedback law in the unconstrained case~\cite{Kalman1960}. 
The number of affine pieces, however, is often large even for low-order examples and short horizons. It is therefore not trivial to calculate the piecewise-affine solution explicitly. 
The first algorithms that were proposed for this task, which are still competitive and available in mature software tools~\cite{Herceg2013,Bemporad2004}, 
exploit the geometric structure of the solution~\cite{Bemporad2002,Tondel2003,Baotic2002}. 
A second, younger class of algorithms is based on finding all optimal active sets of the underlying parametric quadratic program~\cite{Gupta2011,Feller2013,Oberdieck2017,Ahmadi2018,Herceg2015}.  
The set of all optimal active sets determines the solution, 
since every affine piece of the solution is defined by a unique optimal active set under mild conditions (see Sect.~\ref{sec:ProblemStatement} and references given there). 
We call approaches of the second class combinatorial quadratic programming. They are sometimes also called implicit enumeration techniques in the literature. 
We note there exist other approaches than the mentioned ones, notably those based on directional derivatives and vertex enumeration~\cite{Patrinos2010,Monnigmann2012}. 

We propose to combine combinatorial quadratic programming with backward dynamic programming. The essential idea is as follows. 
An optimal active set of the constrained LQR problem with horizon $N$ always contains an optimal active set of the same problem with horizon $N-1$~\cite[Prop.~1]{Monnigmann2019}. 
This implies the set of all optimal active sets for horizon $N$ can be created by copying and extending the optimal active sets for horizon $N-1$. 
The combinatorial complexity of this extension step depends on the number of constraints of the additional stage (denoted $q_{\mathcal{U}\mathcal{X}}$, see Sect.~\ref{sec:ProblemStatement}) but not on the combinatorial complexity of all constraints (a number on the order of $N\cdot q_{\mathcal{U}\mathcal{X}}$). 
As a consequence, the computational effort of the existing combinatorial algorithms can be reduced by building up the set of all optimal active sets by iteratively increasing the horizon to the desired value $N$. Moreover, it is easy to detect that the finitely determined solution has been found if such a solution exists, because, loosely speaking, all known active sets are extended by a stage of inactive constraints in this case (see Prop.~\ref{prop:maxSet}).

Section~\ref{sec:ProblemStatement} summarizes some facts about constrained LQR and combinatorial quadratic programming. 
Sections~\ref{sec:approach} and~\ref{sec:Example} present the proposed algorithm and illustrate it with an example, respectively. 
Conclusions and a brief outlook are given in Sect.~\ref{conclusion}.

\subsection*{Notation}

For any $M\in\mathbb{R}^{a\times b}$ and any ordered set  $\mathcal{M}\subseteq\{1,...,a\}$ 
let $M_{\mathcal{M}}\in\mathbb{R}^{\vert\mathcal{M}\vert\times b}$
be the submatrix of $M$ containing all rows indicated by $\mathcal{M}$. 
Let $\oplus$ and $\ominus$ denote the Minkowski addition and Pontryagin difference, respectively. 

\section{Problem statement and preliminaries} \label{sec:ProblemStatement}

Consider a discrete-time time-invariant linear system
\begin{align} \label{eq:System}
x(k+1)=Ax(k)+Bu(k),\, k= 0, 1, \dots
\end{align}
that must respect constraints of the form
\begin{align*}
  u(k)\in\mathcal{U}\subset\mathbb{R}^{m}, \,
  x(k)\in\mathcal{X}\subset\mathbb{R}^{n}, \,
  k= 0, 1, \dots
\end{align*}
with input variables $u(k)\in\mathbb{R}^m$, state variables $x(k)\in\mathbb{R}^n$, matrices 
$A\in\mathbb{R}^{n\times n}$ and $B\in\mathbb{R}^{n\times m}$, where $(A,B)$ is stabilizable and $\mathcal{U}$ and $\mathcal{X}$ are compact full-dimensional polytopes that contain the origin in their interiors. 

The optimal control problem (OCP) treated in the present paper reads
{\small \begin{align} \label{eq:OCP}
\begin{split}
\min_{U,X} \quad &x(N)^T Px(N)+\sum_{k=0}^{N-1}\left(x(k)^TQx(k)+u(k)^TRu(k)\right) \\
\textrm{s.t.} \quad 
&x(k+1)=Ax(k)+Bu(k), \: k=0,...,N-1\\
&x(k)\in\mathcal{X}, \: k=0,...,N-1\\
&u(k)\in\mathcal{U}, \: k=0,...,N-1\\
&x(N)\in\mathcal{T},
\end{split}
\end{align}}
where $x(0)$ is given, $U=\left(u^T(0),...,u^T(N-1)\right)^T\in\mathbb{R}^{Nm}$ and $X=\left(x^T(1),...,x^T(N)\right)^T\in\mathbb{R}^{Nn}$ collect the inputs and states, respectively, $Q\in\mathbb{R}^{n\times n}$, $Q\succeq 0$ and  $R\in\mathbb{R}^{m\times m}$, $R\succ 0$ are the usual weighting matrices and $N\in\mathbb{N}$ is the horizon. 
We choose $P$ and $K$ to be the optimal cost function matrix and optimal feedback matrix, respectively, of the unconstrained infinite-horizon problem, which implies $P\succ 0$. 
$\mathcal{T}$ is chosen to be the largest possible set with the properties 
$\mathcal{T}\subseteq\mathcal{X}$ (state constraint satisfaction), 
$Kx\in\mathcal{U}$ for all $x\in\mathcal{T}$ (input constraint satisfaction)
and $(A+BK)x\in\mathcal{T}$ for all $x\in\mathcal{T}$ (positive invariance under $K$). 
Let $\mathcal{F}_N$ refer to all $x(0)\in\mathbb{R}^n$ such that~\eqref{eq:OCP} has a solution. Since $\mathcal{F}_N\supseteq\mathcal{T}$ and $\mathcal{T}\ne\emptyset$, $\mathcal{F}_N$ is not empty.

We assume the order of the constraints 
\begin{align} \label{eq:order}
\left.\begin{array}{c}
  u(0)\in\mathcal{U},\quad x(0)\in\mathcal{X},
  \\ 
  \vdots
  \\
  u(N-1)\in\mathcal{U},\quad x(N-1)\in\mathcal{X},
  \\ 
  x(N)\in\mathcal{T}  
\end{array}\right.
\end{align}
and call this the stagewise order with $N+1$ stages $k= 0, \dots, N$. 
Let $q_\mathcal{UX}$ and $q_\mathcal{T}$ be the number of input and state constraints, and terminal constraints, respectively. The total number of constraints $q$ then is $Nq_\mathcal{UX}+q_\mathcal{T}$ for horizon $N$. 
Furthermore, let $\mathcal{Q}=\{1,...,q\}$ and $\mathcal{Q}_0=\{1,...,q_\mathcal{UX}\}$.%
\footnote{The first point in time and thus the first stage is $k=0$, while the first constraint is constraint $i= 1$.  
    Since the initial condition $x(0)$ is usually associated with time $k= 0$, 
		and since the first row of a matrix like $G$ in~\eqref{eq:mpc} is usually considered to be row $i= 1$,
		we have to live with this nuisance. \label{footnote:AnnoyingIndexShift}%
} 

By substituting the dynamics, 
the OCP~\eqref{eq:OCP} can be transformed into a quadratic program (QP) of the form
\begin{align} \label{eq:mpc}
\begin{split}
\min_{U} \quad &\frac{1}{2}x(0)^TYx(0)+x(0)^TFU+\frac{1}{2}U^THU\\
\textrm{s.t.} \quad &GU\leq Ex(0)+w,
\end{split}
\end{align} 
with $Y\in\mathbb{R}^{n\times n}$, $F\in\mathbb{R}^{n\times Nm}$, $H\in\mathbb{R}^{Nm\times Nm}$, $G\in\mathbb{R}^{q\times Nm}$, $w\in\mathbb{R}^{q}$ and $E\in\mathbb{R}^{q\times n}$, 
and where the assumptions on~\eqref{eq:OCP} imply $H\succ 0$~\cite{Bemporad2002}. 
We assume the constraint order from \eqref{eq:order} is preserved in the constraint matrices in~\eqref{eq:mpc}. 

%
For any $x(0)\in\mathcal{F}_N$, 
let $\mathcal{A}(x(0))$ and $\mathcal{I}(x(0))$ refer to the optimal active set $\mathcal{A}(x(0))= \{i\in\mathcal{Q} |G_{\{i\}}U=w_{\{i\}}+E_{\{i\}} x(0)\}$ 
and the corresponding inactive set $\mathcal{I}(x(0))=\mathcal{Q}\backslash\mathcal{A}(x(0))$.
We say $\mathcal{A}(x(0))$ is an \textit{optimal} active set to distinguish it more clearly from \textit{candidate} active sets introduced in Sect.~\ref{sec:combinatorialmpqp}. 

When solving~\eqref{eq:mpc} as a parametric program with parameter $x(0)$, the corresponding optimal control law is a continuous piecewise affine function on a partition of $\mathcal{F}_N$ into full-dimensional polytopes \cite[Sect. 4.1]{Bemporad2002}. 
We denote the set of all optimal active sets $\mathcal{A}$ such that $G_\mathcal{A}$ has full row rank and such $\mathcal{A}$ defines a full-dimensional polytope by $\mathcal{M}_N$.
%
%
We need to consider, however, active sets that define lower dimensional polytopes (such as facets and vertices).
We anticipate all required optimal active sets will be collected in $\mathcal{S}_N\supseteq\mathcal{M}_N$,
which is introduced in Sect.~\ref{sec:approach}.

\subsection{Combinatorial quadratic programming} \label{sec:combinatorialmpqp}

Let $\mathcal{P}(\mathcal{Q})$ refer to the power set of $\mathcal{Q}$ 
and note that the set of active sets $\mathcal{M}_N$ that define the solution is a subset of $\mathcal{P}(\mathcal{Q})$.
It is the basic idea of combinatorial quadratic programming to efficiently determine those $\mathcal{A}\in\mathcal{P}(\mathcal{Q})$ that make up $\mathcal{M}_N$. 

All $\mathcal{A}\in\mathcal{M}_N$ are optimal active sets by definition of $\mathcal{M}_N$. For clarity, 
we call any $\mathcal{A}\in\mathcal{P}(\mathcal{Q})$ that is not known to be optimal a \textit{candidate} active set.  
%
%
A candidate $\mathcal{A}\in\mathcal{P}(\mathcal{Q})$ is optimal, if the linear program (LP)
\begin{subequations} \label{eq:FeasibilityLpWithStationarity}
\begin{align}
\min_{U,x(0),\lambda_{\mathcal{A}},s_{\mathcal{I}},t} \quad &-t \\
\textrm{s.t.} \quad & F^Tx(0) +HU+(G_{\mathcal{A}})^T\lambda_{\mathcal{A}}=0, \label{eq:LPStationarity1} \\
&te_2\leq \lambda_{\mathcal{A}}, \label{eq:LPStationarity2} \\ 
&G_{\mathcal{A}}U-E_{\mathcal{A}}x(0)-w_{\mathcal{A}}=0, \\
&G_{\mathcal{I}}U-E_{\mathcal{I}}x(0)-w_{\mathcal{I}}+s_{\mathcal{I}}=0, \\ 
&te_1\leq s_{\mathcal{I}}, \\
&t\geq 0,
\end{align}
\end{subequations}
has a solution~\cite[Sect.\ 3.1]{Gupta2011}, 
where $e_i=(1\cdots 1)^T$, $i= 1, 2$ are column vectors of appropriate sizes,  
$\lambda_{\mathcal{A}}$ are Lagrangian multipliers and $s_{\mathcal{I}}$ are slack variables. 
Furthermore, we follow~\cite{Gupta2011} in calling an $\mathcal{A}\in\mathcal{P}(\mathcal{Q})$ \textit{feasible} (resp.\ \textit{infeasible}), 
if \eqref{eq:FeasibilityLpWithStationarity} without \eqref{eq:LPStationarity1} and \eqref{eq:LPStationarity2} 
has a solution (resp.\ has no solution). 
Since this smaller LP 
involves a subset of the constraints of~\eqref{eq:FeasibilityLpWithStationarity}, 
feasibility of $\mathcal{A}$ is a necessary condition for optimality of $\mathcal{A}$, or equivalently, an $\mathcal{A}$ that is infeasible is not optimal.  
The LP without \eqref{eq:LPStationarity1} and \eqref{eq:LPStationarity2} is particularly useful, because it typically permits to disregard many candidates $\mathcal{A}\in\mathcal{P}(\mathcal{Q})$ after solving only one LP. 
Specifically, if $\mathcal{A}\in\mathcal{P}(\mathcal{Q})$ is infeasible, 
then additional active constraints cannot result in feasibility and therefore every $\mathcal{A}^\prime\supset\mathcal{A}$ is also infeasible \cite[Thm.~1]{Gupta2011}. 

An optimal active set $\mathcal{A}$ defines a full-dimensional polytope, if $G_\mathcal{A}$ is of full row rank and if both $G_\mathcal{I}U<w_\mathcal{I}+E_\mathcal{I}x(0)$ and $\lambda_\mathcal{A}>0$ hold \cite[Thm. 2]{Tondel2003}. 
Full-dimensional polytopes defined by active sets such that $G_\mathcal{A}$ does not have full row rank are not required, because their polytopes are covered 
by the polytopes defined by the active sets in $\mathcal{M}_N$~\cite[Sect. 3]{Ahmadi2018}.

Algorithm~\ref{algorithm:Gupta} from~\cite{Gupta2011} serves as a reference to which the 
approach proposed here is compared. 
In contrast to the new algorithm stated in Sect.~\ref{sec:approach}, Alg.~\ref{algorithm:Gupta} only considers
candidates $\mathcal{A}$ such that $G_{\mathcal{A}}$ has full rank.  
Since $G_{\mathcal{A}}\in\mathbb{R}^{|\mathcal{A}|\times mN}$, the row rank of $G_{\mathcal{A}}$ is bounded from above by $\min(|A|, mN)$. 
Consequently, $G_{\mathcal{A}}$ does not have full row rank if $|A|>mN$ and only the candidates 
$\mathcal{A}\in\mathcal{P}^\prime(\mathcal{Q})=\{\mathcal{A}\in\mathcal{P}(\mathcal{Q}) | |\mathcal{A}|\le mN\}$ need to be considered (line 2 in Alg.~\ref{algorithm:Gupta}).
The solution $t=0$ to~\eqref{eq:FeasibilityLpWithStationarity} indicates that either $G_{\{i\}}U=w_{\{i\}}+E_{\{i\}}x(0)$ holds for an $i\in\mathcal{I}$ or $\lambda_{\{j\}}=0$ holds for a $j\in\mathcal{A}$. In this case the active set is only added to $\mathcal{M}_N$ if the polytope defined by $\mathcal{A}$ (see e.g.~\cite[Lem. 2]{Jost2015a}) is full-dimensional.

\begin{algorithm2e}[h]
\textbf{Initialization:} set $\mathcal{M}_N=\emptyset$, $\mathcal{S}_N^{\rm pruned}= \emptyset$\\
\For{every $\mathcal{A}_i\in\mathcal{P}'(\mathcal{Q})$ by incr. cardinality}{
	\If{$\mathcal{A}_i\not\supseteq\tilde{\mathcal{A}}$ for all $\tilde{\mathcal{A}}\in\mathcal{S}_N^{\rm pruned}$ and $\text{rowrank}(G_{\mathcal{A}_i})=\vert\mathcal{A}_i\vert$}{
		solve \eqref{eq:FeasibilityLpWithStationarity}\\
		\If{solution $t>0$}{
			add $\mathcal{A}_i$ to $\mathcal{M}_N$
		}
		\ElseIf{solution $t=0$}{
			\If{polytope def.\ by $\mathcal{A}_i$ full-dim.}{
				add $\mathcal{A}_i$ to $\mathcal{M}_N$
			}			
		}
		\Else{
			solve \eqref{eq:FeasibilityLpWithStationarity} without \eqref{eq:LPStationarity1} and \eqref{eq:LPStationarity2}\\
			\If{no solution exists}{
				add $\mathcal{A}_i$ to $\mathcal{S}_N^{\rm pruned}$
			}
		}
	}
}
\textbf{Output:} $\mathcal{M}_N$
\caption{Combinatorial quadratic programming for~\eqref{eq:mpc} \cite{Gupta2011}
\label{algorithm:Gupta}}
\end{algorithm2e}

\section{Dynamic programming approach}\label{sec:approach}

We present an algorithm that combines combinatorial quadratic programming and dynamic programming.  
Essentially, the optimal active sets for horizon $N+1$ can be created by copying and by extending the optimal active sets for horizon $N$.
The computational effort in this extension step is dominated by the combinatorics of the $q_\mathcal{UX}$ constraints of the additional stage,
in contrast to the combinatorics of the $(N+1)q_\mathcal{QX}+ q_\mathcal{T}$ of the total number of constraints for horizon $N+1$.
In addition, it is easy to detect $\mathcal{F}_\infty= \lim_{k\rightarrow\infty} F_k$ has been reached for a finite horizon, i.e., $\mathcal{F}_\infty= \mathcal{F}_N$, 
with the proposed algorithm, 
if such a finite horizon $N$ exists. 

We state two basic relations of active sets for horizons $N$ and $N+1$ in Lems.~\ref{lemma:ShiftAndAugment} and~\ref{lemma:InsertZeroesBeforeT}.
Let $\mathcal{S}_{N}$ refer to the set of all optimal active sets for horizon $N$
which obviously is a superset of $\mathcal{M}_N$.
\begin{lem}[{\cite[Prop. 1]{Monnigmann2019}}]\label{lemma:ShiftAndAugment}
	Consider an OCP \eqref{eq:OCP} and assume its constraints are ordered as in \eqref{eq:order}. 
	Then, for every active set $\mathcal{A}_i\in\mathcal{S}_{N+1}$, there exists an active set $\mathcal{A}_l\in\mathcal{S}_{N}$ such that 
	\begin{align}\label{eq:ShiftAndAugment}
	  \mathcal{A}_i=\mathcal{A}_j\cup (\mathcal{A}_l\oplus\{q_\mathcal{UX}\})
	\end{align}
	for some $\mathcal{A}_j\in\mathcal{P}(\mathcal{Q}_0)$.
\end{lem}

It is evident from Lem.~\ref{lemma:ShiftAndAugment} that combining all active sets in $\mathcal{S}_N$ with all combinations for a single stage $\mathcal{A}_j\in\mathcal{P}(\mathcal{Q}_0)$ results in a superset of all optimal active sets for horizon $N+1$. More formally, $\mathcal{S}_{N+1}\subseteq$
\begin{align}\label{eq:TooLargeSuperset}
  \left\{
    \mathcal{A}_j\cup(\mathcal{A}_l\oplus \{q_\mathcal{UX}\}) |
    \mathcal{A}_j\in\mathcal{P}(\mathcal{Q}_0), \, \mathcal{A}_l\in\mathcal{S}_N
  \right\} .
\end{align}
Some of the sets in~\eqref{eq:TooLargeSuperset} that are not optimal can be removed without having to solve an LP~\eqref{eq:FeasibilityLpWithStationarity}. 
This is stated more precisely in Cor.~\ref{cor:reducedCandidates} below. We need Lem.~\ref{lemma:InsertZeroesBeforeT} as a preparation.
\begin{lem}[{\cite[Lem.~4]{Monnigmann2019}}]\label{lemma:InsertZeroesBeforeT}
	Consider an OCP \eqref{eq:OCP} and assume its constraints are ordered as in \eqref{eq:order}.
	Let $N$ be an arbitrary horizon, let $l\in\mathbb{N}$ be arbitrary, and let $\mathcal{A}_i$ be a candidate active set for horizon $N$ with no active terminal constraints, i.e.,
	\begin{align} \label{eq:InsertZeroesBeforeT}
	  \mathcal{A}_i\subseteq\{1,...,Nq_\mathcal{UX}\}.
	\end{align}  
	Then $\mathcal{A}_i\in\mathcal{S}_N$ if and only if
	$\mathcal{A}_i\in\mathcal{S}_{N+l}$. 
\end{lem}

The set $\mathcal{S}_{N+1}$ can now be constructed from $S_N$ with Lems.~\ref{lemma:ShiftAndAugment} and~\ref{lemma:InsertZeroesBeforeT} in two steps:
(i) copying all $\mathcal{A}\in\mathcal{S}_N$ with no active terminal constraints (Lem.~\ref{lemma:InsertZeroesBeforeT}), 
and (ii) by shifting and augmenting those $\mathcal{A}\in\mathcal{S}_N$ 
that result in active constraints for the terminal stage $N+1$ or the respective previous stage $N$ (Lem.~\ref{lemma:ShiftAndAugment}, specifically~\eqref{eq:ShiftAndAugment}).  
While step (i) always yields optimal active sets for $N+1$, step (ii)
results in candidate active sets. Their optimality still needs to be tested with~\eqref{eq:FeasibilityLpWithStationarity}. 
The sets constructed in steps (i) and (ii) are more formally described in the following corollary.
\begin{cor} \label{cor:reducedCandidates}
	Consider an OCP \eqref{eq:OCP} and assume its constraints are ordered as in \eqref{eq:order}. 
	Assume we know $\mathcal{S}_{N}$. 
	Then 
	\begin{align}\label{eq:reducedCandidatesHelper100}
	  \mathcal{S}_{{N}+1}
	  = 
	  \mathcal{R}^{(1)}
	  \cup
	  \mathcal{R}^{(2)}, 
	\end{align}
	with
	\begin{subequations}
	\begin{align}
    	\mathcal{R}^{(1)}
    	&=
    	\{\mathcal{A}\in\mathcal{S}_{N}\vert\mathcal{A}\subseteq\{1,...,Nq_\mathcal{UX}\}\},
    	\label{eq:reducedCandidatesHelper3}
    	\\
    	\mathcal{R}^{(2)}
		&\subseteq
		\left\{
	    	\mathcal{A}_j\cup(\mathcal{A}_l\oplus \{q_\mathcal{UX}\}) |
	    	\mathcal{A}_j\in\mathcal{P}(\mathcal{Q}_0), \, \mathcal{A}_l\in\tilde{\mathcal{S}}_N
	  	\right\},
	  	\label{eq:reducedCandidates}
	\end{align}
	\end{subequations}
	where $\tilde{\mathcal{S}}_{N}$ contains all elements of $\mathcal{S}_{N}$ that have at least one active constraint in stage $k=N-1$ or $k=N$,
	i.e.,\,
	\begin{align} \label{eq:TildeSN}
		\tilde{\mathcal{S}}_{N}
	  	&=
	  	\{\mathcal{A}\in\mathcal{S}_{N}\vert\mathcal{A}\not\subseteq\{1,...,(N-1)q_\mathcal{UX}\}\}.
	\end{align}
	
	\begin{proof}
		Recall the QP~\eqref{eq:mpc} with horizon $N+1$ has $(N+1)q_{\mathcal{U}\mathcal{X}}+q_\mathcal{T}$ constraints 
		and comprises $N+2$ stages $k= 0, \dots, N, N+1$, where $k= N+1$ corresponds to the terminal constraints		
		(see~\eqref{eq:order} and the subsequent paragraph). 
		The constraints with indices $i>N q_{\mathcal{U}\mathcal{X}}$ belong to the terminal stage (stage $N+1$) 
		and the last stage before the terminal stage (stage $N$).%
		\footnote{See footnote on p.~\ref{footnote:AnnoyingIndexShift}.}
		Now let $\mathcal{A}\in\mathcal{S}_{N+1}$ be arbitrary and distinguish the following two cases from one another: Either
		\begin{align*}
		  \text{(i)}\quad\quad \mathcal{A}\subseteq\{1, \dots, Nq_{\mathcal{U}\mathcal{X}}\}, 
		\end{align*}
		i.e., $\mathcal{A}$ has no active constraints in stages $N$ and $N+1$, or
		\begin{align*}
		  \text{(ii)}\quad\quad \mathcal{A}\not\subseteq\{1, \dots, Nq_{\mathcal{U}\mathcal{X}}\}, 
		\end{align*}
		i.e., $\mathcal{A}$ has at least one active constraint in these two stages. 
		In case (i), condition~\eqref{eq:InsertZeroesBeforeT} is fulfilled.
		Since $\mathcal{A}\in\mathcal{S}_{N+1}$ by assumption, Lem.~\ref{lemma:InsertZeroesBeforeT} applies, which yields $\mathcal{A}\in\mathcal{S}_N$. Together $\mathcal{A}\in\mathcal{S}_N$ and (i) imply $\mathcal{A}\in\mathcal{R}^{(1)}$. 
		In case (ii) we need to show $\mathcal{A}\in\mathcal{R}^{(2)}$. For this purpose, we partition $\mathcal{A}$ into
		\begin{align*}
		  \mathcal{A}_j&= \mathcal{A}\cap \{1, \dots, q_{\mathcal{U}\mathcal{X}}\}
		  \\
		  \bar{\mathcal{A}}_l&= \mathcal{A}\cap \{q_{\mathcal{U}\mathcal{X}}+1, \dots, (N+1)q_{\mathcal{U}\mathcal{X}}+q_\mathcal{T}\}
		\end{align*}
		and let 
		\begin{align*}
		  \mathcal{A}_l&= \bar{\mathcal{A}}_l\ominus \{q_{\mathcal{U}\mathcal{X}}\}.
		\end{align*}
		Then
		\begin{align}\label{eq:reducedCandidatesHelper1}
		  \mathcal{A}= \mathcal{A}_j\cup(\mathcal{A}_l\oplus\{q_{\mathcal{U}\mathcal{X}}\}), \,
		  \mathcal{A}_j\in\mathcal{P}(\mathcal{Q}_0)
		\end{align}
		and because $\mathcal{A}$ contains at least one constraint of the stages $N$ and $N+1$ by assumption (ii), 
		$\bar{\mathcal{A}}_l$ contains at least one constraint of the stages $N$ and $N+1$, 
		and $\mathcal{A}_l$ contains at least one constraint of the stages $N-1$ and $N$. 
		The last statement implies  $\mathcal{A}_l\not\subseteq\{1, \dots, (N-1)q_{\mathcal{U}\mathcal{X}}\}$. 
		Since this is the defining condition of $\tilde{\mathcal{S}}_N$ 
		in~\eqref{eq:TildeSN}, we have 
		\begin{align}\label{eq:reducedCandidatesHelper2}
		  \mathcal{A}_l\in\tilde{\mathcal{S}}_N.
		\end{align}
		Together~\eqref{eq:reducedCandidatesHelper1} and~\eqref{eq:reducedCandidatesHelper2} imply $\mathcal{A}\in\mathcal{R}^{(2)}$.
		We so far showed $\mathcal{A}\in\mathcal{R}^{(1)}$ or $\mathcal{A}\in\mathcal{R}^{(2)}$ for arbitrary $\mathcal{A}\in\mathcal{S}_{N+1}$. 
		It remains to show that the equality in~\eqref{eq:reducedCandidatesHelper3} holds. 
		Since all elements of the right hand side of~\eqref{eq:reducedCandidatesHelper3}  fulfill~\eqref{eq:InsertZeroesBeforeT} and $\mathcal{A}\in\mathcal{S}_N$, Lem.~\ref{lemma:InsertZeroesBeforeT} applies and~\eqref{eq:InsertZeroesBeforeT} holds and all these elements are also elements of $\mathcal{S}_{N+1}$, which completes the proof.
	\end{proof}
\end{cor}

Note that active sets $\mathcal{A}\in\mathcal{S}_N$ with no active terminal constraints and at least one active constraint in stage $N-1$ are treated with both~\eqref{eq:reducedCandidatesHelper3} and~\eqref{eq:reducedCandidates}. Nevertheless, $\mathcal{R}^{(1)}\cap\mathcal{R}^{(2)}=\emptyset$ holds, because the concerned active sets are treated differently (copied in~\eqref{eq:reducedCandidatesHelper3} and extended in~\eqref{eq:reducedCandidates}).

If a finite horizon $N$ exists such that the $\mathcal{F}_\infty= \mathcal{F}_N$, it is easy to detect this $N$ has been reached. More precisely, the following statement holds.
\begin{prop}\label{prop:maxSet}
Consider an OCP \eqref{eq:OCP} and assume its constraints are ordered as in \eqref{eq:order}. Assume we know $\mathcal{S}_N$ and let $\tilde{\mathcal{S}}_N$ be defined as in Cor.~\ref{cor:reducedCandidates}. 

If $\tilde{\mathcal{S}}_{N}=\emptyset$, then the solution for the finite horizon $N$ 
defined by $\mathcal{S}_N$ is the solution for all horizons $\tilde{N}\ge N$. 
Furthermore, the corresponding optimal control law is identical for all horizons $\tilde{N}\ge N$.

	\begin{proof}
		First note that $\tilde{\mathcal{S}}_{N}=\emptyset$ implies that the right hand side of \eqref{eq:reducedCandidates} and therefore $\mathcal{R}^{(2)}$ is empty. It follows that
		\begin{align}\label{eq:FiniteFHelper1}
		\mathcal{S}_{{N}+1} =\mathcal{R}^{(1)}
		\end{align}
		with~\eqref{eq:reducedCandidatesHelper100}. 
		Secondly, by definition of $\tilde{\mathcal{S}}_N$ in Cor.~\ref{cor:reducedCandidates}, $\tilde{\mathcal{S}}_{N}=\emptyset$ implies there exist no active sets in $\mathcal{S}_N$ with active terminal constraints, which yields
		\begin{align}\label{eq:FiniteFHelper2}
		\mathcal{R}^{(1)}=\mathcal{S}_N.
		\end{align}
		Combining~\eqref{eq:FiniteFHelper1} and~\eqref{eq:FiniteFHelper2} yields $\mathcal{S}_{N+1}=\mathcal{S}_N$.
		$\mathcal{S}_{N+l}= \mathcal{S}_N$ for all $l\in\mathbb{N}$ follows by induction. 
		Finally, an active set $\mathcal{A}\in\mathcal{S}_N$ with no active terminal constraints not only reappears in $\mathcal{S}_{N+l}$ for all $l\in\mathbb{N}$
		but such a $\mathcal{A}$ is known to define the same polytope and optimal affine feedback law for all horizons $N+l$~\cite[Prop.~6]{Monnigmann2019}.
	\end{proof}
\end{prop}


Corollary~\ref{cor:reducedCandidates} and Prop.~\ref{prop:maxSet} are illustrated with an example in Sect.~\ref{sec:illustration}.

\subsection{Implementational aspects}\label{subsec:ImplementationalAspects}

Assuming the set $S_N$ is known for some horizon $N$, the set $S_{N+1}$ can be determined with Alg.~\ref{algorithm:SNp1FromSN}, 
which essentially implements Cor.~\ref{cor:reducedCandidates}. 
Specifically, lines 4,5 correspond to~\eqref{eq:reducedCandidatesHelper3} and lines 8,9 correspond to~\eqref{eq:reducedCandidates}
in Cor.~\ref{cor:reducedCandidates}. 
\begin{algorithm2e}[h]
\SetKwInOut{Initialization}{Initialization}
\SetKwInOut{Input}{Input}
\SetKwInOut{Output}{Output}
\textbf{Input:} $\mathcal{S}_N$, $\mathcal{S}_N^{\rm degen.}$\\
\textbf{Initialization:} set $\mathcal{S}_{N+1}=\emptyset$, $\mathcal{S}_{N+1}^{\rm degen.}$ and $\mathcal{S}_{N+1}^{\rm pruned}= \emptyset$\\
\For{every $\mathcal{A}_l\in\mathcal{S}_{N}$}{
	\If{$\mathcal{A}_l\subseteq\{1,...,Nq_\mathcal{UX}\}$}{
		add $\mathcal{A}_l$ to $\mathcal{S}_{N+1}$\\
		\If{$\mathcal{A}_l\in\mathcal{S}_N^{\rm degen.}$}{
			add $\mathcal{A}_l$ to $\mathcal{S}_{N+1}^{\rm degen.}$
		}
	}
	\If{$\mathcal{A}_l\not\subseteq\{1,...,(N-1)q_\mathcal{UX}\}$}{
		\For{every $\mathcal{A}_i=\mathcal{A}_j\cup(\mathcal{A}_l\oplus\{q_\mathcal{UX}\})$ with $\mathcal{A}_j\in\mathcal{P}(\mathcal{Q}_0)$  by increasing cardinality}{
			\If{$\mathcal{A}_i\not\supseteq\tilde{\mathcal{A}}$ for all $\tilde{\mathcal{A}}\in\mathcal{S}_{N+1}^{\rm pruned}$}{
				solve \eqref{eq:FeasibilityLpWithStationarity} for QP for horizon $N+1$\\
				\If{solution exists}{
					add $\mathcal{A}_i$ to $\mathcal{S}_{N+1}$\\
					\If{solution $t=0$}{
						add $\mathcal{A}_i$ to $\mathcal{S}_{N+1}^{\rm degen.}$
					}
				}
				\Else{
					solve \eqref{eq:FeasibilityLpWithStationarity} without \eqref{eq:LPStationarity1} and \eqref{eq:LPStationarity2} for QP for horizon $N+1$\\
					\If{no solution exists}{
						add $\mathcal{A}_i$ to $\mathcal{S}_{N+1}^{\rm pruned}$ 
					}
				}
			}
		}
	}
}
\textbf{Output:} $\mathcal{S}_{N+1}$, $\mathcal{S}_{N+1}^{\rm degen.}$
\caption{Determination of $\mathcal{S}_{N+1}$ from $\mathcal{S}_N$
\label{algorithm:SNp1FromSN}}
\end{algorithm2e}
Candidate active sets are tested for optimality for the OCP with horizon $N+1$,
unless they can be disregarded because they are supersets of a known infeasible active set (line 10). 
Candidate active sets are added to $\mathcal{S}_{N+1}$ if they are optimal (lines 12,13) and tested for feasibility otherwise (line 17).
Sets known to be infeasible are stored (lines 18,19) in order to be able to quickly dismiss supersets that appear later. 
The set $\mathcal{S}_N^{\rm degen.}$ collects all $\mathcal{A}\in\mathcal{S}_N$ such that the solution to~\eqref{eq:FeasibilityLpWithStationarity} is $t=0$. $\mathcal{S}_{N+1}^{\rm degen.}$ is then determined by collecting all $\mathcal{A}\in\mathcal{S}_N^{\rm degen.}$ that were copied with~\eqref{eq:reducedCandidatesHelper3} (lines 6,7) and all candidate active sets such that the solution to~\eqref{eq:FeasibilityLpWithStationarity} is $t=0$ (lines 14,15).
Since active sets that are infeasible for horizon $N$ are not necessarily infeasible for horizon $N+1$, 
$\mathcal{S}_{N+1}^{\rm pruned}$ is initialized with the empty set in line 2. 

The initial sets $\mathcal{S}_1$ and $\mathcal{S}_1^{\rm degen.}$ can be determined with Alg.~\ref{algorithm:initHorizon},
which proceeds analogously to Alg.~\ref{algorithm:Gupta} but does not discard candidate active sets such that $G_\mathcal{A}$ is not of full rank. In case the solution to~\eqref{eq:FeasibilityLpWithStationarity} is $t=0$, the candidate active set is added to both $\mathcal{S}_1$ and $\mathcal{S}_1^{\rm degen.}$.
\begin{algorithm2e}[h]
\textbf{Initialization:} set $\mathcal{S}_1=\emptyset$, $\mathcal{S}_1^{\rm degen.}=\emptyset$ and $S_1^{\rm pruned}= \emptyset$\\
\For{every $\mathcal{A}_i\in\mathcal{P}(\{1,...,q_\mathcal{UX}+q_\mathcal{T}\})$ by incr. cardinality}{
	\If{$\mathcal{A}_i\not\supseteq\tilde{\mathcal{A}}$ for all $\tilde{\mathcal{A}}\in\mathcal{S}_1^{\rm pruned}$}{
		solve \eqref{eq:FeasibilityLpWithStationarity} for QP with horizon $1$\\
		\If{solution exists}{
			add $\mathcal{A}_i$ to $\mathcal{S}_1^{\rm degen.}$\\
			\If{solution $t=0$}{
				add $\mathcal{A}_i$ to $\mathcal{S}_1^{\rm degen.}$
			}
		}
		\Else{
			solve \eqref{eq:FeasibilityLpWithStationarity} without \eqref{eq:LPStationarity1} and \eqref{eq:LPStationarity2} for QP with horizon $1$\\
			\If{no solution exists}{
				add $\mathcal{A}_i$ to $S_1^{\rm pruned}$
			}
		}
	}
}
\textbf{Output:} $\mathcal{S}_1$, $\mathcal{S}_1^{\rm degen.}$
\caption{Determination of $\mathcal{S}_1$
\label{algorithm:initHorizon}}
\end{algorithm2e}

The overall dynamic programming approach is stated in Alg.~\ref{algorithm:dynamicProgrammingApproach}.
The algorithm terminates if the desired horizon $N_{\rm max}$ has been reached 
or if an $N$ such that $\mathcal{F}_\infty= \mathcal{F}_N$ has been found with Prop.~\ref{prop:maxSet}. 
The condition in line 5 of Alg.~\ref{algorithm:dynamicProgrammingApproach} merely is a compact way of stating that all active sets in $\mathcal{S}_{N+1}$ 
have no active constraints in stages $N$ and $N+1$ and thus is equivalent to $\tilde{\mathcal{S}}_{N+1}=\emptyset$.
Lines 8-14 in Alg.~\ref{algorithm:dynamicProgrammingApproach} reduce $\mathcal{S}_N$ to $\mathcal{M}_N\subseteq\mathcal{S}_N$ by discarding all active sets such that $G_\mathcal{A}$ does not have full rank and by testing all active sets that are element of $\mathcal{S}_N^{\rm degen.}$ for defining a full-dimensional polytope.
\begin{algorithm2e}[h]
\textbf{Input:} $\mathcal{S}_1$, $\mathcal{S}_1^{\rm degen.}$ (from Alg.~\ref{algorithm:initHorizon}), $N_{\max}\ge 1$\\
\textbf{Initialization:} set $\mathcal{M}_N=\emptyset$\\
\For{$N=1$ to $N_{\max}-1$}{
	determine $\mathcal{S}_{N+1}$ and $\mathcal{S}_{N+1}^{\rm degen.}$ with Alg.~\ref{algorithm:SNp1FromSN}\\
	\If{$\mathcal{S}_{N+1}\subseteq\mathcal{P}(\{1,...,Nq_\mathcal{UX}\})$}{
		break 
	}
}
$N=N+1$\\
\For{every $\mathcal{A}_k\in\mathcal{S}_N$}{
	\If{$\text{rowrank}(G_{\mathcal{A}_k})=\vert\mathcal{A}_k\vert$}{
		\If{$\mathcal{A}_k\in\mathcal{S}_N^{\rm degen.}$}{
			\If{polytope def.\ by $\mathcal{A}_k$ full-dim.}{
				add $\mathcal{A}_k$ to $\mathcal{M}_N$
			}
		}
		\Else{
			add $\mathcal{A}_k$ to $\mathcal{M}_N$
		}
	}
}
\textbf{Output:} $\mathcal{M}_N$
\caption{Dynamic programming approach to solving constrained linear-quadratic OCPs
\label{algorithm:dynamicProgrammingApproach}}
\end{algorithm2e}

\begin{figure}[thpb]
   \begin{center}
   \begin{minipage}[t]{0.43\textwidth}
   \includegraphics[width=\textwidth]{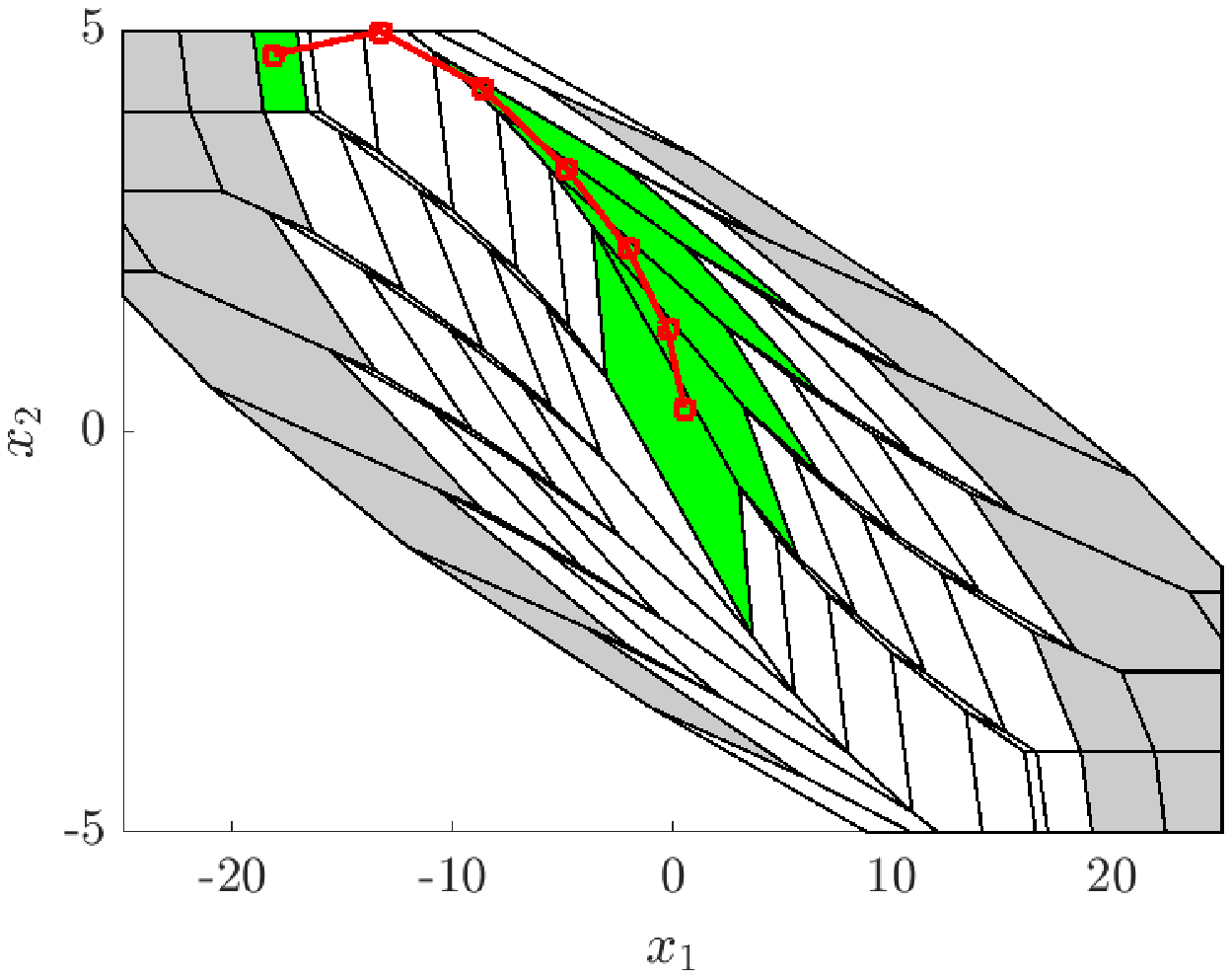}
   \end{minipage}\\
   \begin{minipage}[t]{0.43\textwidth}
   \includegraphics[width=\textwidth]{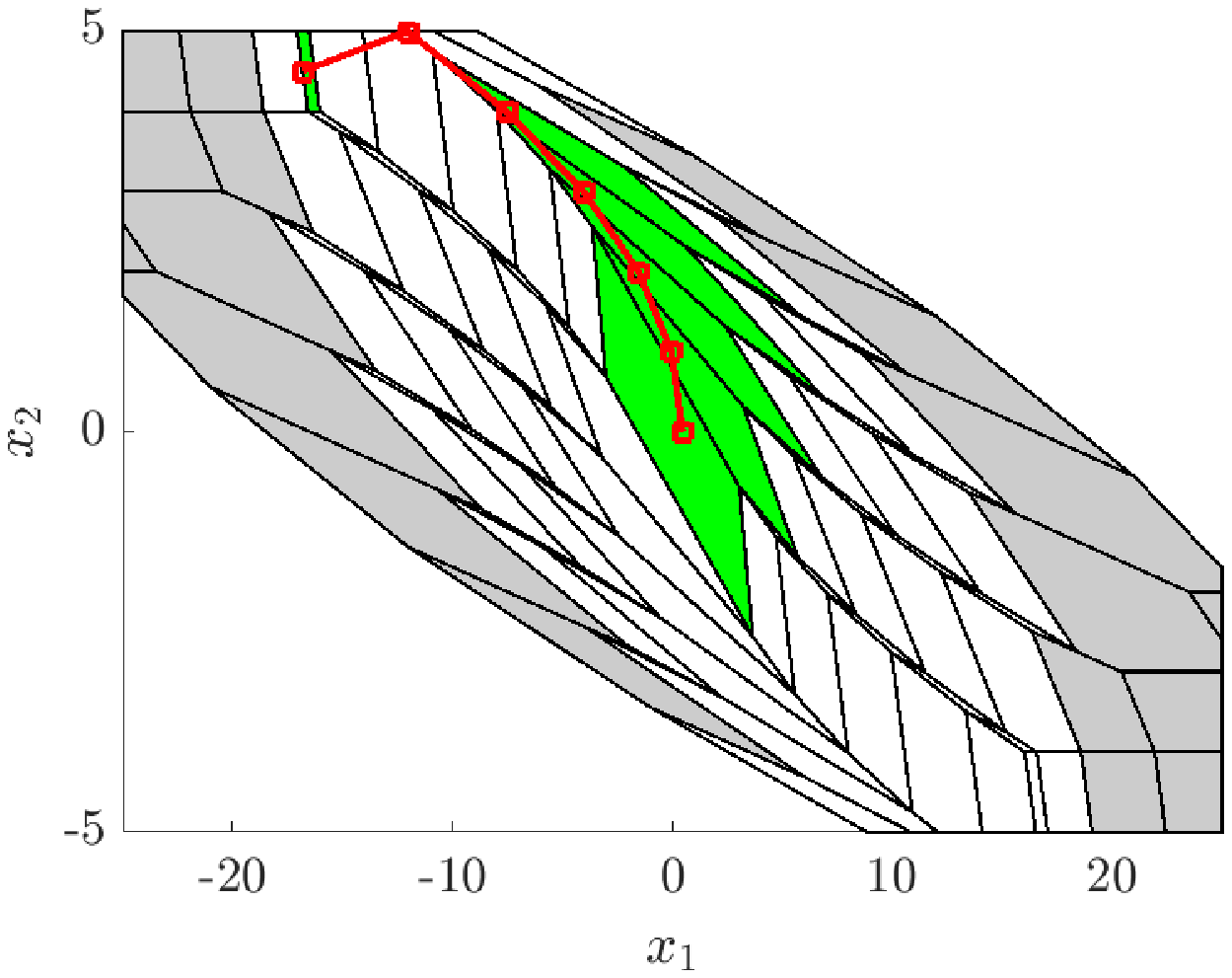}
   \end{minipage}
   \caption{Solution for OCP~\eqref{eq:OCP} with the example from Sect.~\ref{sec:Example}, $N=6$. Red lines and green polytopes mark sample closed-loop trajectories and the polytopes they pass through, respectively. Polytopes with active terminal constraints are gray.
   Top: $\mathcal{A}_i=\{12,13,19,25,31\}\in\mathcal{S}_6$ where $G_{\mathcal{A}i}$ is of full rank defines a full-dimensional polytope
(leftmost green polytope). 
   $\mathcal{A}_i\in\mathcal{S}_6$ results with~\eqref{eq:ShiftAndAugment} from $\mathcal{A}_l=\{6,7,13,19,25\}\in\mathcal{S}_{5}$.
   $\mathcal{A}_l\in\mathcal{S}_5$ does not exceed the maximal cardinality for horizon $N= 5$ but
    $G_{\mathcal{A}l}$ is not of full rank.
   The corresponding polytope is a 1-dimensional facet. 
   Bottom: The active set $\mathcal{A}_i=\{7,12,13,19,25,31\}$ where $G_{\mathcal{A}i}$ is of full rank defines a full-dimensional polytope
    (leftmost green polytope). $\mathcal{A}_i\in\mathcal{S}_N$ is constructed with~\eqref{eq:ShiftAndAugment} from the active set $\mathcal{A}_l=\{1,6,7,13,19,25\}\in\mathcal{S}_5$. $\mathcal{A}_l$ exceeds the maximal cardinality for horizon $N=5$. Therefore $G_{\mathcal{A}l}$ is not of full rank and $\mathcal{A}_l$ defines
    a 1-dimensional facet.}
   \label{fig:illustrationLICQfail}
   \end{center}
\end{figure}

In contrast to Alg.~\ref{algorithm:Gupta}, 
Alg.~\ref{algorithm:dynamicProgrammingApproach} does not consider all $\mathcal{A}\in\mathcal{P}'(\mathcal{Q})$, 
but it generates candidate active sets with Cor.~\ref{cor:reducedCandidates}, thus reducing their number.  
On the other hand, candidate active sets such that $G_\mathcal{A}$ is not of full rank are discarded 
in Alg.~\ref{algorithm:Gupta}, but not in Alg.~\ref{algorithm:dynamicProgrammingApproach}. 
In fact, an $\mathcal{A}_i\in\mathcal{S}_{N+1}$ such that $G_{\mathcal{A}i}$ is of full rank
may result with~\eqref{eq:ShiftAndAugment}
from a $\mathcal{A}_l\in\mathcal{S}_N$ such that $G_{\mathcal{A}l}$ is not of full rank
This is illustrated in Fig.~\ref{fig:illustrationLICQfail}.
The results in Sect.~\ref{sec:omputationalEffort} show that considerable overall savings result with the method proposed here even though active sets such that $G_\mathcal{A}$ is not of full rank are no longer discarded.

\section{Example} \label{sec:Example}

We illustrate Cor.~\ref{cor:reducedCandidates} and Prop.~\ref{prop:maxSet} in Sect.~\ref{sec:illustration} and subsequently 
analyze the computational effort of the new approach in Sect.~\ref{sec:omputationalEffort}.
We use the double integrator~\cite{Gutman1987} 
\begin{align*}
x(k+1)=\left(\begin{array}{cc}1&1\\0&1\end{array}\right) x(k)+\left(\begin{array}{c}0.5\\1\end{array}\right)u(k)
\end{align*}
with input constraints $\vert u(k)\vert \leq 1$, state constraints $\vert x_1(k)\vert \leq 25$, $\vert x_2(k)\vert \leq 5$ and cost function matrices 
$Q=1\in\mathbb{R}^{2\times 2}$, $R= 0.1$ as an example.
The terminal cost $P$ and set $\mathcal{T}$ are as described in Sect.~\ref{sec:ProblemStatement}. 

\subsection{Illustration of Cor.~\ref{cor:reducedCandidates} and Prop.~\ref{prop:maxSet}}
\label{sec:illustration}
\begin{figure}[thpb]
   \begin{center}
\begin{minipage}[t]{0.43\textwidth}
\includegraphics[trim=0px 0 0px 0, clip, width=\textwidth]{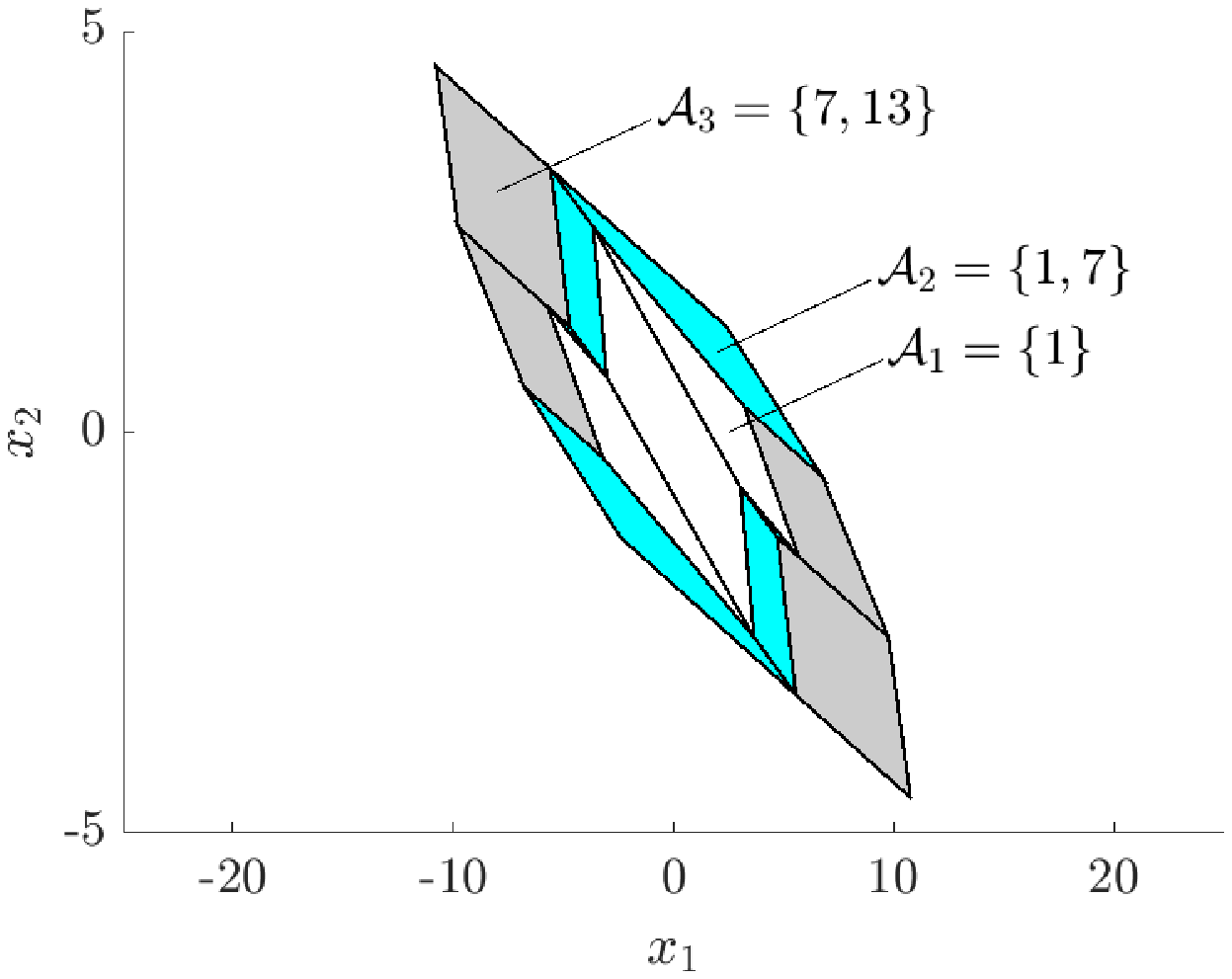}
\end{minipage}
\begin{minipage}[t]{0.43\textwidth}
\includegraphics[trim=0px 0 0px 0, clip, width=\textwidth]{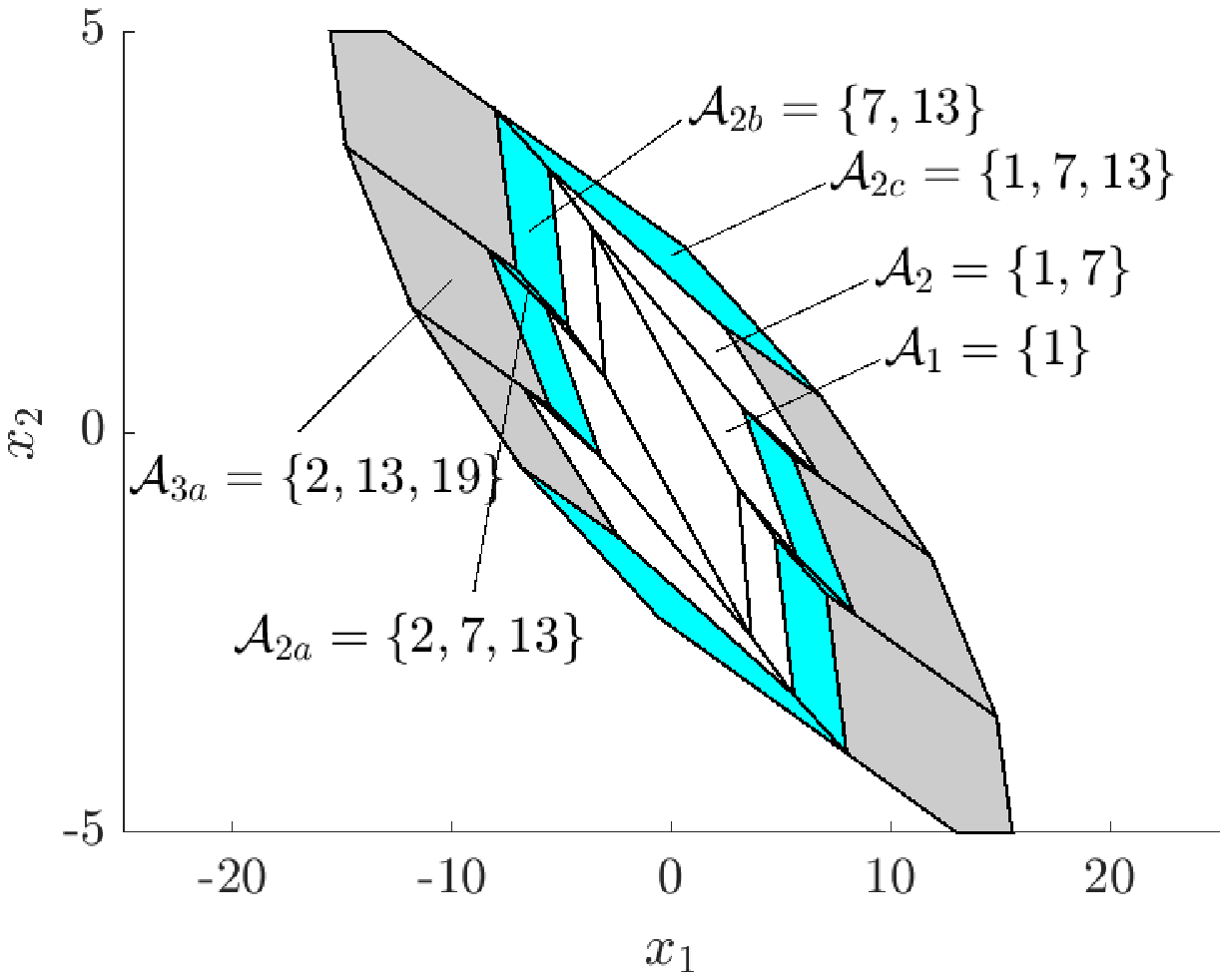}
\end{minipage}
   \caption{Solution for OCP~\eqref{eq:OCP} for $N=2$ (top) and $N=3$ (bottom). 
     Gray polytopes have at least one active terminal constraint, 
     blue polytopes have no active terminal constraints but at least one active constraints in stage $N-1$,
     white polytopes have no active constraints in stages $N$ and $N-1$.}
     \label{fig:buildSolution}
     \end{center}
\end{figure}

\begin{figure}[thpb]
   \begin{center}
\begin{minipage}[t]{0.43\textwidth}
\includegraphics[trim=0px 0 0px 0, clip, width=\textwidth]{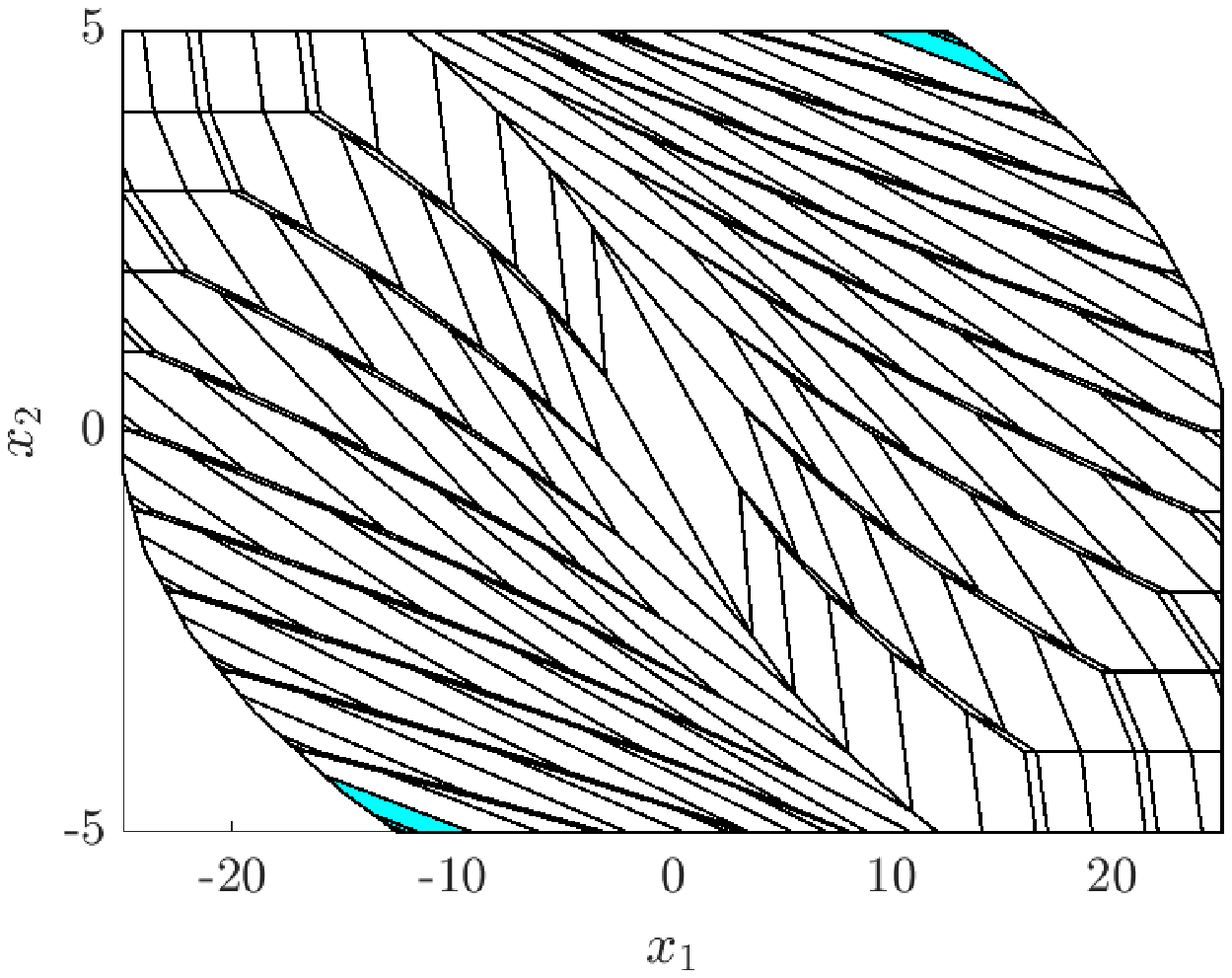}
\end{minipage}
\begin{minipage}[t]{0.43\textwidth}
\includegraphics[trim=0px 0 0px 0, clip, width=\textwidth]{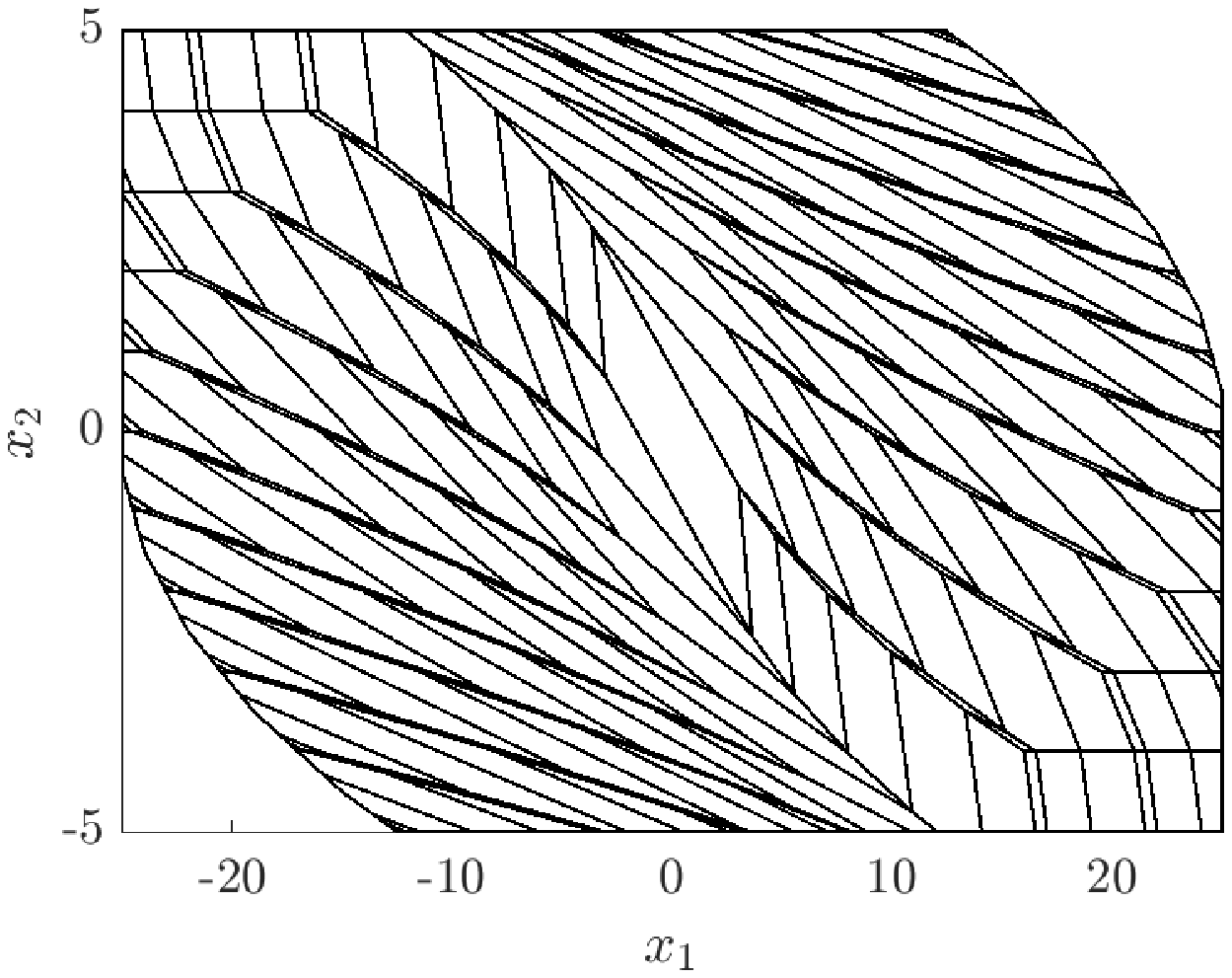}
\end{minipage}
   \caption{Solution for OCP~\eqref{eq:OCP} for $N=15$ (top) and $N=16$ (bottom). 
   Colors are as in Fig.~\ref{fig:buildSolution}.}
      \label{fig:buildSolutionInfty}
   \end{center}
\end{figure}

Figures~\ref{fig:buildSolution} and \ref{fig:buildSolutionInfty} show solutions for the OCP~\eqref{eq:OCP} as a function of the horizon.
Gray polytopes correspond to active sets with at least one active terminal constraint.
Blue polytopes correspond to active sets with no active terminal constraints (stage $N$) but at least one active constraint in stage $N-1$.
White polytopes correspond to active sets with no active constraints in stages $N-1$ and $N$.
The set $\mathcal{S}_N$ contains the active sets of all shown polytopes.
$\tilde{\mathcal{S}}_N$ as defined in~\eqref{eq:TildeSN} contains the active sets of all blue and gray polytopes. 

Figure~\ref{fig:buildSolution} illustrates Cor. \ref{cor:reducedCandidates} with the solutions for $N= 2$ and $N= 3$.
All elements in $\mathcal{S}_2$ with no active terminal constraints (blue and white polytopes) also appear in the solution for the increased horizon $\mathcal{S}_3$.
For example, $\mathcal{A}_1$ and $\mathcal{A}_2$ exists for both $N= 2$ and $N= 3$ according to Lem.~\ref{lemma:InsertZeroesBeforeT} and~\eqref{eq:reducedCandidatesHelper3}. 
All elements in $\tilde{\mathcal{S}}_2$ (blue and gray polytopes) are the basis for the generation of candidate active sets for $\mathcal{S}_{3}$ according to~\eqref{eq:reducedCandidates} in Cor.~\ref{cor:reducedCandidates}. For example, 
the blue polytopes with active sets $\mathcal{A}_{2a},...,\mathcal{A}_{2c}$ in $\mathcal{S}_3$ result with~\eqref{eq:reducedCandidates} from 
$\mathcal{A}_2$, and the gray polytope $\mathcal{A}_{3a}$ results with~\eqref{eq:reducedCandidates} 
from $\mathcal{A}_3$. 
 
Proposition~\ref{prop:maxSet} states the solution will not change for horizons $\tilde{N}> N$, if $\tilde{\mathcal{S}}_N=\emptyset$. The horizon $N= 16$ is the shortest horizon such that $\tilde{\mathcal{S}}_N= \emptyset$. The resulting solution set $\mathcal{F}_\infty= \mathcal{F}_{16}$ is shown in Fig.~\ref{fig:buildSolutionInfty}. All polytopes belong to active sets with no active constraints in stages $N-1$ and $N$ (i.e., white polytopes) as expected. We note that $\mathcal{F}_\infty= \mathcal{F}_{15}$ holds, since all active sets in $\mathcal{M}_{15}$, where $\mathcal{M}_{15}$ results from Alg.~\ref{algorithm:dynamicProgrammingApproach}, have no active terminal constraints~\cite[Prop.~6]{Monnigmann2019}.

\subsection{Computational effort} 
\label{sec:omputationalEffort}

\begin{figure}[tbh]
   \begin{center}
   \begin{subfigure}[b]{0.45\textwidth}
        \begin{center}
		\includegraphics[trim=0px 0 0px 0, clip, width=\textwidth]{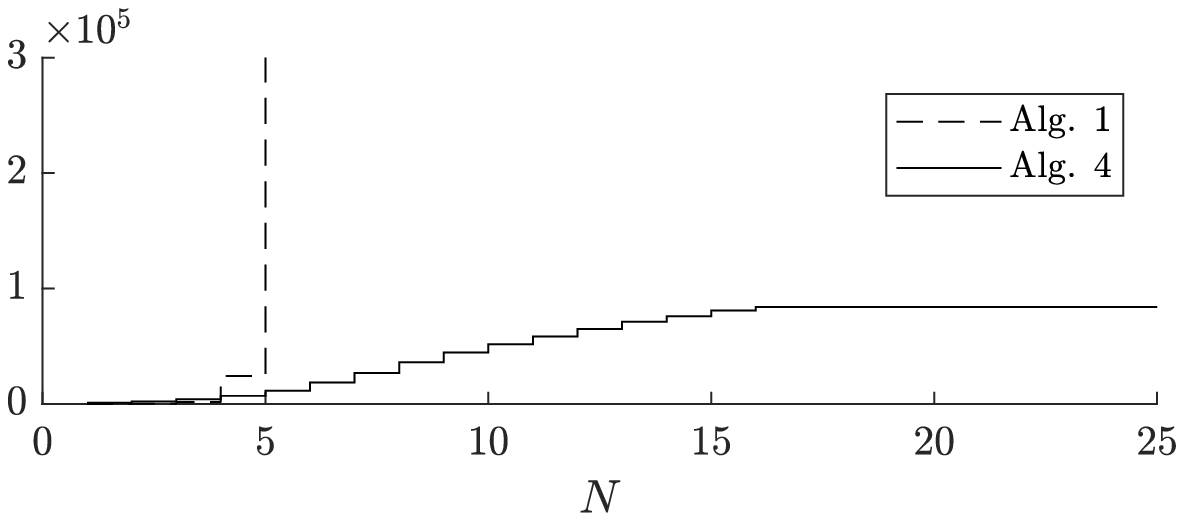}
        \caption{candidate active sets}
        \label{fig:Example1Candidates}
        \end{center}
	\end{subfigure}\\
	\vspace{2mm}
   \begin{subfigure}[b]{0.23\textwidth}
        \begin{center}
		\includegraphics[trim=0px 0 0px 0, clip, width=\textwidth]{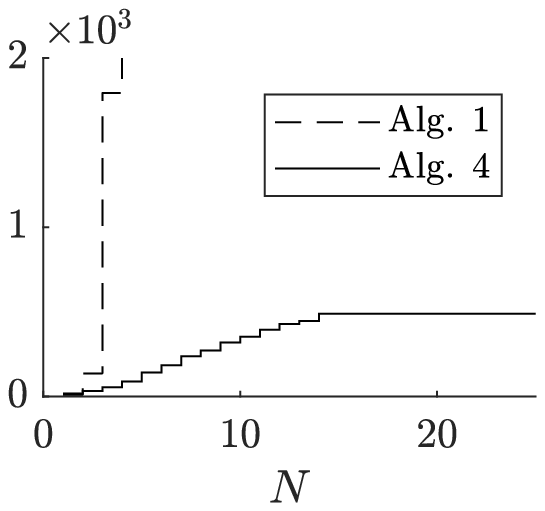}
        \caption{rank tests}
        \label{fig:Example1LICQ}
        \end{center}
	\end{subfigure}
   \begin{subfigure}[b]{0.23\textwidth}
        \begin{center}
		\includegraphics[trim=0px 0 0px 0, clip, width=\textwidth]{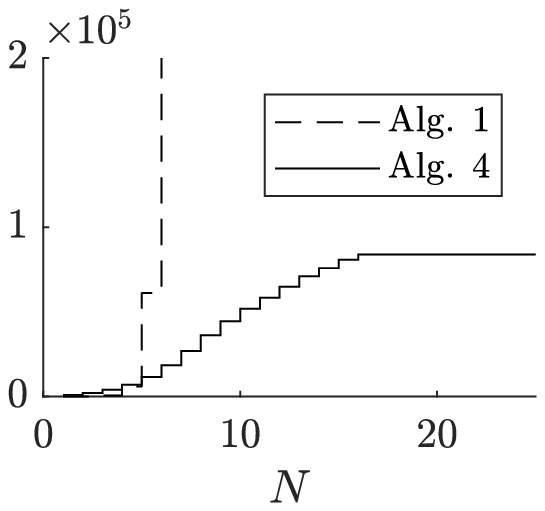}
        \caption{pruning tests}
        \label{fig:Example1Pruned}
        \end{center}
	\end{subfigure}\\
	\vspace{2mm}
   \begin{subfigure}[b]{0.23\textwidth}
        \begin{center}
		\includegraphics[trim=0px 0 0px 0, clip, width=\textwidth]{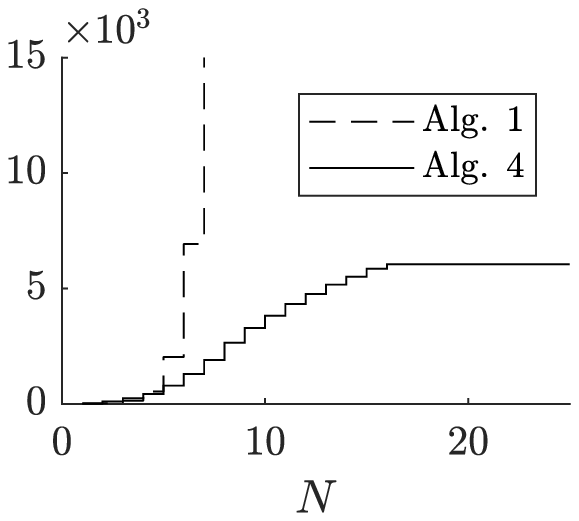}
        \caption{optimality tests}
        \label{fig:Example1Optimality}
        \end{center}
	\end{subfigure}
   \begin{subfigure}[b]{0.23\textwidth}
        \begin{center}
		\includegraphics[trim=0px 0 0px 0, clip, width=\textwidth]{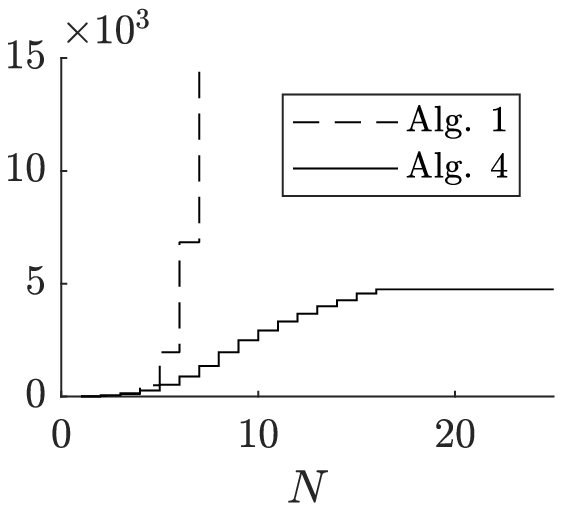}
        \caption{feasibility tests}
        \label{fig:Example1Feasibility}
        \end{center}
	\end{subfigure}
    \caption{Number of candidate active sets, rank tests, pruning tests, optimality tests, feasibility tests for Algs.~\ref{algorithm:Gupta} and~\ref{algorithm:dynamicProgrammingApproach}.}
    \label{fig:Example1}
    \end{center}
\end{figure}

Figure~\ref{fig:Example1} compares the computational costs of the existing approach (Alg.~\ref{algorithm:Gupta}) and the approach proposed here (Alg.~\ref{algorithm:dynamicProgrammingApproach}). 
The figure shows the numbers of generated candidate active sets, pruning tests, 
rank tests, optimality tests with \eqref{eq:FeasibilityLpWithStationarity}, and feasibility tests with \eqref{eq:FeasibilityLpWithStationarity} without \eqref{eq:LPStationarity1} and \eqref{eq:LPStationarity2} for both algorithms.
In Alg.~\ref{algorithm:Gupta}, the number of generated candidate active sets equals the cardinality of $\mathcal{P}^\prime(\mathcal{Q})$.
In Alg.~\ref{algorithm:dynamicProgrammingApproach}, the number of generated candidate active sets is equal to the sum of the cardinalities of $\mathcal{P}(\{1,...,q_\mathcal{UX}+q_\mathcal{T}\})$ and all $\mathcal{R}^{(2)}$ for horizons $N=1,...,N_{\max}-1$. 

\begin{figure}[tbh]
   \begin{center}
   \includegraphics[trim=0px 0 0px 0, clip, width=0.45\textwidth]{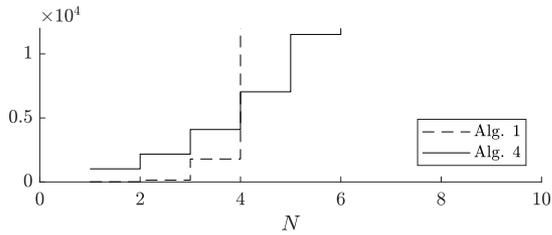}
   \caption{Number of candidate active sets for Algs.~\ref{algorithm:Gupta} and~\ref{algorithm:dynamicProgrammingApproach} (detail of Fig.~\ref{fig:Example1Candidates}).}
   \label{fig:Example1CandidatesDetail}
   \end{center}
\end{figure}

All curves in Fig.~\ref{fig:Example1} intersect for small values of $N$, which implies Alg.~\ref{algorithm:Gupta} outperforms 
Alg.~\ref{algorithm:dynamicProgrammingApproach} for small $N$.  
The point of intersection is enlarged in Fig.~\ref{fig:Example1CandidatesDetail}.   
For larger $N$, Alg.~\ref{algorithm:dynamicProgrammingApproach} is more efficient than Alg.~\ref{algorithm:Gupta} and it is evident from Fig.~\ref{fig:Example1} this difference gets more pronounced as $N$ increases. 

Algorithms~\ref{algorithm:Gupta} and~\ref{algorithm:dynamicProgrammingApproach} mainly differ with respect 
to the fundamentally different procedures for generating candidate active sets. This difference becomes evident in Fig.~\ref{fig:Example1Candidates}.
Most importantly, a plateau results for Alg.~\ref{algorithm:dynamicProgrammingApproach} in Fig.~\ref{fig:Example1} because $\mathcal{F}_\infty= \mathcal{F}_{15}$. 
Since fewer candidate active sets entail fewer rank, pruning, optimality and feasibility tests, the qualitative difference 
of the curves in Fig.~\ref{fig:Example1Candidates} is inherited by the remaining curves. 

%
In both algorithms, the optimality tests with LP~\eqref{eq:FeasibilityLpWithStationarity} only need to be carried out for a subset of the candidate active sets (compare Figs.~\ref{fig:Example1Candidates} and~\ref{fig:Example1Optimality}). In Alg.~\ref{algorithm:Gupta}, candidates can be dismissed because they are supersets of known infeasible active sets and because of the rank test (line 3). In the new approach, in contrast, candidates are only dismissed with the first criterion (line 10 in Alg.~\ref{algorithm:SNp1FromSN}). However, even though two criteria can be used in Alg.~\ref{algorithm:Gupta} to reduce the number of optimality tests, a smaller number of optimality tests still results in Alg.~\ref{algorithm:dynamicProgrammingApproach} (Fig.~\ref{fig:Example1Optimality}), 
because a considerably smaller number of candidate active sets must be tested to begin with (Fig.~\ref{fig:Example1Candidates}). 
This effect carries over to the number of feasibility tests (Fig.~\ref{fig:Example1Feasibility}). 

Rank tests appear only in the last step of Alg.~\ref{algorithm:dynamicProgrammingApproach} (line 9). 
Consequently, their overall number is also smaller than in Alg.~\ref{algorithm:Gupta} (Fig.~\ref{fig:Example1LICQ}). 

We note the LPs solved in Alg.~\ref{algorithm:dynamicProgrammingApproach} have fewer constraints than the LPs that are solved in Alg.~\ref{algorithm:Gupta}, because the number of constraints increases with the horizon, which is always $N$ in Alg.~\ref{algorithm:Gupta} while it increases from 1 to $N$ in Alg.~\ref{algorithm:dynamicProgrammingApproach}. 

Finally, it should be noted that Alg.~\ref{algorithm:dynamicProgrammingApproach} requires more memory than Alg.~\ref{algorithm:Gupta}, because the set 
$\mathcal{S}_N\supseteq\mathcal{M}_N$ is required in Alg.~\ref{algorithm:dynamicProgrammingApproach} but not in Alg.~\ref{algorithm:Gupta}.

\section{Conclusions} \label{conclusion}
We introduced a new algorithm for determining the set of optimal active sets that determine the solution to the constrained LQR problem. 
It is the central idea of the proposed algorithm to build active sets by iteratively increasing the horizon of the constrained LQR problem. In doing so, the combinatorial complexity of existing algorithms is greatly reduced. The anticipated reduction was illustrated with an example.

Extending the optimal active sets for horizon $N-1$ to those for $N$ formally corresponds to a backward dynamic programming step~\cite{Monnigmann2019}. 
It is an obvious question to ask whether also the geometric approaches (see Sect.~\ref{sec:Intro} for a brief summary)
could built up the solution by iteratively increasing the horizon. This is hampered by the fact that the optimal feedback law and its state space partition for horizon $N-1$ is not in general contained in the optimal feedback law for horizon $N$~\cite{Munoz2004}. 
Since it is easy to identify the persistent polytopes 
with Lems.~\ref{lemma:InsertZeroesBeforeT} and~\ref{lemma:ShiftAndAugment}, future work will reconsider combining backward dynamic programming with the existing geometric approaches.

\section{Acknowledgements}                           
This work was supported by the German Federal Ministry for Economic Affairs and Energy under grant 0324125C. 

\bibliographystyle{plain}        
\bibliography{root}           

\end{document}